\documentclass{article}
\usepackage{amsmath}
\usepackage{amsfonts}
\usepackage{amssymb}
\usepackage{graphicx}
\newtheorem{theorem}{Theorem}

\newtheorem{corollary}[theorem]{Corollary}

\newtheorem{lemma}[theorem]{Lemma}

\newenvironment{proof}[1][Proof]{\textbf{#1.} }{\ \rule{0.5em}{0.5em}}

\begin{document}

\title{Quantum Spheres As Groupoid C*-algebras}
\author{Albert J. L. Sheu\thanks{Partially supported by NSF-Grant DMS-9303231}\\Department of Mathematics\\University of Kansas\\Lawrence, KS 66045-2142\\U. S. A.}
\date{}
\maketitle
\begin{abstract}
In this paper, we show that the C*-algebras of quantum spheres can be realized
as concrete groupoid C*-algebras.
\end{abstract}

\section{Introduction}

Ever since the successful use of groupoid C*-algebras \cite{Re} by Connes in
studying the interaction between operator algebras and foliation geometry
\cite{Co}, it has been well recognized that the theory of groupoid C*-algebras
provides a very useful tool and a nearly universal context for the study of
operator algebras, especially in connection with geometric structures. The
results obtained in \cite{CuM,MRe,SSU,Sh:rd,Sh:cqg} are a few examples
supporting this aspect. In particular, it has been found that some interesting
C*-algebras associated with singular foliations can be effectively described
and studied in the context of C*algebras of (essentially) discrete groupoids
\cite{Sh:sf,Sh:cqg}.

Recent fast developments in deformation quantization
\cite{Ri:dq,Ri:dqar,Sh:qp} including quantum groups
\cite{D,FRT,So,Wo:ts,Wo:cm,Ri:cq} showed that the problem of quantization of
Poisson structures and group structures is intimately related to the
underlying singular symplectic foliation \cite{We,LuWe}. Furthermore in the
paper \cite{Sh:cqg}, it is found that the C*-algebra arising from deformation
quantization using quantum groups can be studied by the method of groupoids.

In fact, it is shown \cite{Sh:cqg} that C*-algebras $C(SU(n)_{q})$ of the
quantum groups $SU(n)_{q}$ can be realized as a C*-subalgebra of some groupoid
C*-algebras. We also know that $C(SU(2)_{q})\simeq C^{*}(\frak{F})$ where
\[
\frak{F}=\{(m,j,k)\;|\;\text{ if }k=\infty\text{, then }m=j\}
\]
is a subgroupoid of the Toeplitz groupoid ${\mathbb{Z}}\times({\mathbb{Z}%
}\times\overline{{\mathbb{Z}}}|_{\overline{{\mathbb{Z}}}_{\geq}})$.

In the study of quantizing Poisson Lie groups, it is also of great interest to
see how the Poisson structures reduce to homogeneous spaces and how the
resulting Poisson spaces can be quantized \cite{DoGu}. Fundamental examples
are provided by the study of quantum spheres \cite{Po,VaSo}. It is also found
that C*-algebras arising from such a quantization are closely related to the
underlying symplectic structure and the corresponding groupoid \cite{Sh:cqg}.
In particular, we know that $C(S_{q}^{2n+1})$ of the quantum sphere
$S_{q}^{2n+1}$ \cite{VaSo} can be realized as a C*-subalgebra of some groupoid
C*-algebra, and $C(S_{\mu c}^{2})\simeq C^{*}(\frak{G})$ where $S_{\mu c}^{2}$
is Podles' quantum sphere \cite{Po} with $c>0$ and
\[
\frak{G}=\{(j,j,k_{1},k_{2})\;|\;\text{ }k_{1}=\infty\text{ or }k_{2}%
=\infty\}
\]
is a subgroupoid of the 2-dimensional Toeplitz groupoid ${\mathbb{Z}}%
^{2}\times\overline{{\mathbb{Z}}}^{2}|_{\overline{{\mathbb{Z}}}_{\geq}^{2}}$.

However it is not clear and it is an interesting question whether other
$C(S_{q}^{2n+1})$\ or $C(SU(n)_{q})$ can be realized as a groupoid C*-algebra
itself instead of a C*-subalgebra of some groupoid C*-algebra. In this paper,
we give a positive answer to this question for $C(S_{q}^{2n+1})$ and show that
the underlying groupoid can be identified explicitly and is independent of the
parameter $q$. This also gives a proof of the fact that for a fixed $n$, the
C*-algebras $C(S_{q}^{2n+1})$ are isomorphic and independent of $q$. We
mention that the above statements are also true for quantum $SU\left(
3\right)  $, but the proof is more complicated and we will deal with it in a
separate paper.

\section{Quantum sphere and groupoid}

We refer to the paper \cite{Sh:cqg} for notations and concepts.

Recall that the C*-algebra of quantum group $SU(n)_{q}$ is generated by
elements $u_{ij}$ satisfying certain commutation relations and the C*-algebra
of quantum spheres $S_{q}^{2n+1}=SU(n)_{q}\backslash SU(n+1)_{q}$ defined as
homogeneous quantum spaces \cite{N} can be identified with \textbf{\ }
\[
C(S_{q}^{2n+1})=C^{*}(\{u_{n+1,m}|\,1\leq m\leq n+1\}).
\]
In \cite{Sh:cqg}, we give explicit embeddings of $C\left(  SU(n)_{q}\right)  $
and $C(S_{q}^{2n+1})$ into concrete groupoid C*-algebras.

Recall that $C(S_{q}^{2n+1})$ can be embedded in the groupoid C*-algebra
$C^{*}(\mathcal{F}^{n})$, where $\mathcal{F}^{n}$ is the $n$-dimensional
Toeplitz groupoid ${\mathbb{Z}}\times({\mathbb{Z}}^{n}\times\overline
{{\mathbb{Z}}}^{n}|_{\overline{{\mathbb{Z}}}_{\geq}^{n}})$. Let
\[
F_{i}:=\{w\in F:w_{i}=\infty\}\cong\overline{{\mathbb{Z}}}_{\geq}^{n-1}
\]
denote a `face' of the unit space $F\cong\overline{{\mathbb{Z}}}_{\geq}^{n}$
of the groupoid $\mathcal{F}^{n}$ and $\partial F=\cup_{i=1}^{n}F_{i}$ be the
`boundary' of $F$.

Let $\rho_{i*}:\mathcal{F}^{n}|_{F_{i}}\rightarrow\mathcal{F}^{n}$ be the
proper continuous groupoid inclusion map and
\[
\rho_{i}=C^{*}\left(  \rho_{i*}\right)  :C^{*}(\mathcal{F}^{n})\rightarrow
C^{*}(\mathcal{F}^{n})|_{F_{i}}=C^{*}(\mathcal{F}^{n}|_{F_{i}})
\]
be the corresponding C*-algebra homomorphism implementing the restriction to
the face $F_{i}$. Similarly we consider the embeddings $\rho_{\partial
F*}:\mathcal{F}^{n}|_{\partial F}\rightarrow\mathcal{F}^{n}$ and
$\rho_{i,\partial F*}:\mathcal{F}^{n}|_{F_{i}}\rightarrow\mathcal{F}%
^{n}|_{\partial F}$, and the corresponding algebra homomorphisms
$\rho_{\partial F}$ and $\rho_{i,\partial F}$.

Now the restriction of $C(S_{q}^{2n+1})$ to the face $F_{n}$ dominates the
restriction to any other $F_{i}$, in the sense that $\rho_{i}$ on
$C(S_{q}^{2n+1}) $ factors through $\rho_{n}$ on $C(S_{q}^{2n+1})$. More
precisely, we have the following lemma.%

%TCIMACRO{\TeXButton{BeginLemma}{\begin{lemma}}}%
%BeginExpansion
\begin{lemma}%
%EndExpansion
$\rho_{n,\partial F}$ when restricted to $C(S_{q}^{2n+1})|_{\partial F}$ is a
C*-algebra isomorphism
\[
\rho_{n,\partial F}=C^{*}\left(  \rho_{n,\partial F*}\right)  :C(S_{q}%
^{2n+1})|_{\partial F}\rightarrow C(S_{q}^{2n+1})|_{F_{n}}
\]
with its inverse
\[
\pi_{n,\partial F}=C^{*}\left(  \pi_{n,\partial F*}\right)  :C(S_{q}%
^{2n+1})|_{F_{n}}\rightarrow C(S_{q}^{2n+1})|_{\partial F}
\]
arising from the proper continuous groupoid quotient map $\pi_{n,\partial
F*}:\mathcal{F}^{n}|_{\partial F}\rightarrow\mathcal{F}^{n}|_{F_{n}}$ sending
$\left(  z,x,w\right)  $ to $\left(  z,x,w^{\prime}\right)  $ where
$w_{n}^{\prime}=\infty$ and $w_{j}^{\prime}=w_{j}$ for all $j<n$.%

%TCIMACRO{\TeXButton{EndLemma}{\end{lemma}}}%
%BeginExpansion
\end{lemma}%
%EndExpansion%

%TCIMACRO{\TeXButton{Proof}{\proof}}%
%BeginExpansion
\proof
%EndExpansion

In fact, since the restriction $\mathcal{A}|_{\partial F}$ of any
C*-subalgebra $\mathcal{A}\subset C^{*}(\mathcal{F}^{n})$ can be identified,
under the map $\oplus_{i=1}^{n}\rho_{i,\partial F}:\mathcal{A}|_{\partial
F}\rightarrow\oplus_{i=1}^{n}\mathcal{A}|_{F_{i}}$, with the pull-back
\[
\mathcal{A}^{\prime}=\left\{  \oplus_{i}a_{i}\in\oplus_{i}\mathcal{A}|_{F_{i}%
}:\rho_{ij}\left(  a_{j}\right)  =\rho_{ji}\left(  a_{i}\right)  \text{ for
all }i,j\right\}
\]
of the homomorphisms $\mathcal{A}|_{F_{j}}\overset{\rho_{ij}}{\rightarrow
}\mathcal{A}|_{F_{i}\cap F_{j}}$ arising from the proper continuous embedding
$\rho_{ij*}$ of the face $F_{i}\cap F_{j}$ into the face $F_{i}$, we only need
to show that in the commuting diagram
\[%
\begin{array}
[c]{rcl}%
\mathcal{A}|_{\partial F} & \overset{\rho_{i,\partial F}}{\rightarrow} &
\mathcal{A}|_{F_{i}}\\
\rho_{n,\partial F}\downarrow &  & \downarrow\rho_{ni}\\
\mathcal{A}|_{F_{n}} & \overset{\rho_{in}}{\rightarrow} & \mathcal{A}%
|_{F_{i}\cap F_{n}}%
\end{array}
\]
with $\mathcal{A}=C(S_{q}^{2n+1})$, $\rho_{ni}$ is an isomorphism for $i<n$
with the inverse
\[
\pi_{ni}=C^{*}\left(  \pi_{ni*}\right)  :C(S_{q}^{2n+1})|_{F_{i}\cap F_{n}%
}\rightarrow C(S_{q}^{2n+1})|_{F_{i}}
\]
arising from the groupoid quotient map $\pi_{ni*}:\mathcal{F}^{n}|_{F_{i}%
}\rightarrow\mathcal{F}^{n}|_{F_{i}\cap F_{n}}$ sending $\left(  z,x,w\right)
$ to $\left(  z,x,w^{\prime}\right)  $ where $w_{j}^{\prime}=w_{j}$ for all
$j<n$ and $w_{n}^{\prime}=\infty$. Note that $\pi_{ni}\circ\rho_{in}%
=C^{*}\left(  \pi_{n,\partial F*}|_{F_{i}}\right)  $.

It is clear that $\rho_{ni}\circ\pi_{ni}=Id$ since $\pi_{ni\ast}\circ
\rho_{ni\ast}=Id$. On the other hand, we can verify $\pi_{ni}\circ\rho
_{ni}=Id$ on the generators $\rho_{i}\left(  u_{n+1,m}\right)  $ of
$C(S_{q}^{2n+1})|_{F_{i}}$. Recall that $u_{n+1,m}\in C(S_{q}^{2n+1})$ is
identified with
\[
(\tau_{n+1}\otimes\pi_{L_{1}})(\Delta u_{n+1,m})=t_{n+1}\otimes\underset
{n+1-m}{\underbrace{\gamma\otimes...\otimes\gamma}}\otimes\underset
{m-1}{\underbrace{\alpha^{\ast}\otimes1\otimes...\otimes1}}
\]
in $C^{\ast}(\mathcal{F}^{n})$. It is easy to see that $\rho_{i}((\tau
_{n+1}\otimes\pi_{L_{1}})\Delta(u_{n+1,m}))\ =0$ for $m\leq n+1-i$, and hence
we only need to show that $\left(  \pi_{ni}\circ\rho_{ni}\right)  \left(
\rho_{i}(u_{n+1,m})\ \right)  =\rho_{i}(u_{n+1,m})\ $ for $m>n+1-i$. But
\[
\rho_{i}\left(  (\tau_{n+1}\otimes\pi_{L_{1}})(\Delta u_{n+1,n+2-i})\right)
=t_{n+1}\otimes\underset{i-1}{\underbrace{\gamma\otimes...\otimes\gamma}%
}\otimes\underset{m-1}{\underbrace{s_{i}^{-1}\otimes1\otimes...\otimes1}}
\]
with $s_{i}$ the canonical generator of $C\left(  \mathbb{T}\right)  \cong
C^{\ast}\left(  \left\{  \infty\right\}  \times\mathbb{Z}\right)  $, and
\[
\rho_{i}\left(  (\tau_{n+1}\otimes\pi_{L_{1}})(\Delta u_{n+1,m})\right)
=t_{n+1}\otimes\underset{n+1-m}{\underbrace{\gamma\otimes...\otimes\gamma}%
}\otimes\underset{m-1}{\underbrace{\alpha^{\ast}\otimes1\otimes...\otimes1}}
\]
for $m>n+2-i$, are of the form $a\otimes1$ where $a\in C^{\ast}\left(
\mathcal{F}^{n-1}|_{F_{i}}\right)  $ and $1\in C^{\ast}\left(  \mathbb{Z}%
\times\overline{\mathbb{Z}}|_{\overline{\mathbb{Z}}_{\geq}}\right)  $. Since
clearly $\left(  \pi_{ni}\circ\rho_{ni}\right)  \left(  a\otimes1\right)
=a\otimes1$, we get $\left(  \pi_{ni}\circ\rho_{ni}\right)  \left(  \rho
_{i}(u_{n+1,m})\ \right)  =\rho_{i}(u_{n+1,m})\ $ for $m>n+1-i$, and hence
$Id=\pi_{ni}\circ\rho_{ni}$ on $\rho_{i}\left(  C(S_{q}^{2n+1})\right)
=C(S_{q}^{2n+1})|_{F_{i}}$ for $i<n$. (Note that $\pi_{n_{i}}$ and
$\pi_{n,\partial F}$ are well-defined because $C\left(  S_{q}^{2n+1}\right)
|_{F_{i}}\subset C^{\ast}\left(  \mathcal{F}^{n-1}|_{F_{i}}\right)  \otimes1$
for $i<n$.)%

%TCIMACRO{\TeXButton{End Proof}{\endproof}}%
%BeginExpansion
\endproof
%EndExpansion

\section{The groupoid realization of quantum sphere}

By the fact \cite{Sh:cqg} that $C(\mathbb{T})\otimes\mathcal{K}\left(
\ell^{2}\left(  \mathbb{Z}_{\geq}^{n}\right)  \right)  \cong C^{\ast}\left(
\mathcal{F}^{n}|_{\mathbb{Z}_{\geq}^{n}}\right)  \subset C(S_{q}^{2n+1})$, we
have the commuting diagram
\[%
\begin{array}
[c]{ccccccccc}%
0 & \rightarrow &  C(\mathbb{T})\otimes\mathcal{K}\left(  \ell^{2}\left(
\mathbb{Z}_{\geq}^{n}\right)  \right)  & \rightarrow &  C(S_{q}^{2n+1}) &
\overset{\rho_{\partial F}}{\rightarrow} & C(S_{q}^{2n+1})|_{\partial F} &
\rightarrow & 0\\
&  & \parallel &  & \cap &  & \cap &  & \\
0 & \rightarrow &  C(\mathbb{T})\otimes\mathcal{K}\left(  \ell^{2}\left(
\mathbb{Z}_{\geq}^{n}\right)  \right)  & \rightarrow &  C^{\ast}%
(\mathcal{F}^{n}) & \overset{\rho_{\partial F}}{\rightarrow} & C^{\ast
}(\mathcal{F}^{n}|_{\partial F}) & \rightarrow & 0
\end{array}
\]
with exact rows. It is clear that the structure of $C(S_{q}^{2n+1})$ is
completely determined by $C(S_{q}^{2n+1})|_{\partial F}\subset C^{\ast
}(\mathcal{F}^{n}|_{\partial F})$. On the other hand, $C(S_{q}^{2n+1}%
)|_{\partial F}=\pi_{n,\partial F}\left(  C(S_{q}^{2n+1})|_{F_{n}}\right)  $
can be determined explicitly when $C(S_{q}^{2n+1})|_{F_{n}}$ is.

Define a subquotient groupoid $\frak{F}_{n}$ of $\mathcal{F}^{n}={\mathbb{Z}%
}\times({\mathbb{Z}}^{n}\times\overline{{\mathbb{Z}}}^{n}|_{\overline
{{\mathbb{Z}}}_{\geq}^{n}})$ as follows. Let%

\[
\widetilde{\frak{F}_{n}}:=\{\left(  z,x,w\right)  \in\mathcal{F}^{n}%
|\;w_{i}=\infty\Longrightarrow
\]
\[
x_{i}=-z-x_{1}-x_{2}-...-x_{i-1}\text{ and }x_{i+1}=...=x_{n}=0\}
\]
be a subgroupoid of $\mathcal{F}^{n}$. Define $\frak{F}_{n}:=\widetilde
{\frak{F}_{n}}/\sim$ where $\sim$ is the equivalence relation generated by
\[
(z,x,w)\sim(z,x,w_{1},...,w_{i}=\infty,\infty,...,\infty)
\]
for all $(z,x,w)$ with $w_{i}=\infty$ for an $1\leq i\leq n$.%

%TCIMACRO{\TeXButton{BeginTheorem}{\begin{theorem}}}%
%BeginExpansion
\begin{theorem}%
%EndExpansion
$C(S_{q}^{2n+1})\simeq C^{*}(\frak{F}_{n})$.%

%TCIMACRO{\TeXButton{EndTheorem}{\end{theorem}}}%
%BeginExpansion
\end{theorem}%
%EndExpansion%

%TCIMACRO{\TeXButton{Proof}{\proof}}%
%BeginExpansion
\proof
%EndExpansion

Actually we prove by induction on $n$ that
\[
C(S_{q}^{2n+1})=\iota_{n}\eta_{n}\left(  C^{\ast}(\frak{F}_{n})\right)  \simeq
C^{\ast}(\frak{F}_{n})
\]
where $\eta_{n}$ is the injective map arising from the groupoid quotient map
$\eta_{n\ast}:\widetilde{\frak{F}_{n}}\rightarrow\frak{F}_{n}$ and $\iota
_{n}:C^{\ast}\left(  \widetilde{\frak{F}_{n}}\right)  \rightarrow C^{\ast
}\left(  \mathcal{F}^{n}\right)  $ implements the inclusion of the open
subgroupoid $\widetilde{\frak{F}_{n}}$ in $\mathcal{F}^{n}$.

When $n=1$, the equivalence relation $\sim$ is trivial, i.e. $\widetilde
{\frak{F}_{1}}=\widetilde{\frak{F}_{1}}/\sim=\frak{F}_{1}$ and hence
$\iota_{1}\eta_{1}\left(  C^{*}(\frak{F}_{1})\right)  =$ $\iota_{1}\left(
C^{*}(\widetilde{\frak{F}_{1}})\right)  $ with
\[
\widetilde{\frak{F}_{1}}=\{(m,j,k)\;|\;\text{ if }k=\infty\text{, then
}m=j\}.
\]
On the other hand, the C*-algebra $C(S_{q}^{2n+1})=C\left(  SU\left(
2\right)  _{q}\right)  \subset C^{*}(\mathcal{F}^{1})$ is generated by
\[
(\tau_{2}\otimes\pi_{L_{1}})(u_{2,1})=t_{2}\otimes\gamma
\]
and
\[
(\tau_{2}\otimes\pi_{L_{1}})(u_{2,2})=t_{2}\otimes\alpha^{*},
\]
and hence coincides with $\iota_{1}\left(  C^{*}(\widetilde{\frak{F}_{1}%
})\right)  $ as can be easily verified.

Now we assume the induction hypothesis
\[
C(S_{q}^{2n-1})=\iota_{n-1}\eta_{n-1}\left(  C^{*}(\frak{F}_{n-1})\right)
\simeq C^{*}(\frak{F}_{n-1}).
\]

By identifying $\left(  z,x,w\right)  \in\mathcal{F}^{n-1}\subset
\mathbb{Z}\times\mathbb{Z}^{n-1}\times\overline{\mathbb{Z}}_{\geq}^{n-1}$ with
$\left(  z,x,0,w,\infty\right)  \in\mathcal{F}^{n}|_{F_{n}}\subset
\mathbb{Z}\times\mathbb{Z}^{n}\times\overline{\mathbb{Z}}_{\geq}^{n-1}%
\times\left\{  \infty\right\}  $, we regard $\mathcal{F}^{n-1}$ and hence
$\widetilde{\frak{F}_{n-1}}$ as open subgroupoids of $\mathcal{F}^{n}|_{F_{n}%
}$, and regard $C^{*}\left(  \mathcal{F}^{n-1}\right)  $ and $C^{*}\left(
\widetilde{\frak{F}_{n-1}}\right)  $ as C*-subalgebras of $C^{*}\left(
\mathcal{F}^{n}\right)  $. Note that under this identification,
\[
\widetilde{\frak{F}_{n-1}}:=\{(z,x,w)\in\mathcal{F}^{n}|\;w_{n}=\infty\text{,
}x_{n}=0\text{, and for }i<n\text{, }
\]
\[
\text{``}w_{i}=\infty\Longrightarrow x_{i}=-z-x_{1}-...-x_{i-1}\text{ and
}x_{i+1}=...=x_{n}=0\text{''}\}
\]
and $\frak{F}_{n-1}=\widetilde{\frak{F}_{n-1}}/\sim$ with respect to the same
equivalence relation $\sim$.

Now note that $C(S_{q}^{2n+1})|_{F_{n}}\cong C(S_{q}^{2n-1})\subset
C^{*}\left(  \mathcal{F}^{n-1}\right)  $ under the automorphism $\phi$ on
$C^{*}\left(  \mathcal{F}^{n}|_{F_{n}}\right)  \supset C^{*}\left(
\mathcal{F}^{n-1}\right)  $, corresponding to the groupoid isomorphism
$\phi_{*}$ on $\mathcal{F}^{n}|_{F_{n}}\subset\mathbb{Z}\times\mathbb{Z}%
^{n}\times\overline{\mathbb{Z}}_{\geq}^{n-1}$ defined by
\[
\phi_{*}\left(  z,x,w\right)  =\left(  z,x+\left(  0,...,0,-z-\sum_{i=1}%
^{n-1}x_{i}\right)  ,w\right)  .
\]
More precisely, $C(S_{q}^{2n+1})|_{F_{n}}=\phi\left(  C(S_{q}^{2n-1})\right)
$.

Thus from the induction hypothesis, we get
\[
C(S_{q}^{2n+1})|_{F_{n}}=\phi\left[  \iota_{n-1}\eta_{n-1}\left(
C^{*}(\frak{F}_{n-1})\right)  \right]  .
\]

Since under the groupoid isomorphism $\phi_{*}$, $\widetilde{\frak{F}_{n-1}} $
is mapped to the open subgroupoid
\[
\widetilde{\frak{F}_{n-1}^{\prime}}:=\phi_{*}\left(  \widetilde{\frak{F}%
_{n-1}}\right)  =\{(z,x,w)\in\mathcal{F}^{n}|\;w_{n}=\infty\text{, }%
x_{n}=-z-\sum_{i=1}^{n-1}x_{i}\text{, }
\]
\[
\text{and for }i<n\text{, ``}w_{i}=\infty\Longrightarrow x_{i}=-z-x_{1}%
-...-x_{i-1}\text{ and }x_{i+1}=...=x_{n}=0\text{''}\}
\]
\[
=\{(z,x,w)\in\mathcal{F}^{n}|\text{ }w_{n}=\infty\text{, and for any }i\leq
n\text{,}\;
\]
\[
\text{ ``}w_{i}=\infty\Longrightarrow x_{i}=-z-x_{1}-...-x_{i-1}\text{ and
}x_{i+1}=...=x_{n}=0\text{''}\}
\]
of $\mathcal{F}^{n}|_{F_{n}}$ and the equivalence relation $\sim$ is
preserved, it is clear that
\[
C(S_{q}^{2n+1})|_{F_{n}}=\iota_{n-1}^{\prime}\eta_{n-1}^{\prime}\left(
C^{*}(\frak{F}_{n-1}^{\prime})\right)
\]
where $\eta_{n-1*}^{\prime}:\widetilde{\frak{F}_{n-1}^{\prime}}\rightarrow
\frak{F}_{n-1}^{\prime}:=\widetilde{\frak{F}_{n-1}^{\prime}}/\sim$ is the
quotient map of the subgroupoid $\widetilde{\frak{F}_{n-1}^{\prime}}$ of
$\mathcal{F}^{n}|_{F_{n}}$, and $\iota_{n-1}^{\prime}:C^{*}\left(
\widetilde{\frak{F}_{n-1}^{\prime}}\right)  \rightarrow C^{*}\left(
\mathcal{F}^{n}|_{F_{n}}\right)  $ is the inclusion.

Now for $\left(  z,x,w\right)  \in\widetilde{\frak{F}_{n-1}^{\prime}}$, we
have $w_{n}=\infty$, and if $w_{i}=\infty$ for some $i<n$ then $x_{n}=0$.
Recall that $\pi_{n,\partial F\ast}\left(  z,x,w\right)  =\left(
z,x,w^{\prime}\right)  $ with $w^{\prime}=\left(  w_{1},...,w_{n-1}%
,\infty\right)  $. We claim that
\[
\widetilde{\frak{F}_{n-1}^{\prime\prime}}:=\pi_{n,\partial F\ast}^{-1}\left(
\widetilde{\frak{F}_{n-1}^{\prime}}\right)  =\left\{  \left(  z,x,w\right)
\in\mathcal{F}^{n}|_{\partial F}:\text{for any }i\leq n\text{,}\right\}
\]
\[
\text{ ``}w_{i}=\infty\Longrightarrow x_{i}=-z-x_{1}-...-x_{i-1}\text{ and
}x_{i+1}=...=x_{n}=0\text{''}\}.
\]
Indeed, if $\left(  z,x,w\right)  \in\mathcal{F}^{n}|_{\partial F}$ and
$\pi_{n,\partial F\ast}\left(  z,x,w\right)  \in\widetilde{\frak{F}%
_{n-1}^{\prime}}$, then $\left(  z,x,w\right)  $ satisfies the above condition
since it is already satisfied by $\left(  z,x,w^{\prime}\right)  \in
\widetilde{\frak{F}_{n-1}^{\prime}}$ and $w$ does not have more $\infty$
entries than $w^{\prime}$. On the other hand, if $\left(  z,x,w\right)
\in\mathcal{F}^{n}|_{\partial F}$ satisfies the above condition, then either
$\left(  z,x,w\right)  \in\mathcal{F}^{n}|_{F_{n}}$ and hence $\left(
z,x,w^{\prime}\right)  =\left(  z,x,w\right)  $ since $w_{n}=\infty$, or
$w_{i}=\infty$ for some $i<n$ which implies that $x_{i}=-z-x_{1}-...-x_{i-1}$
with $x_{i+1}=...=x_{n}=0$ and hence
\[
x_{n}=0=\left(  -z-x_{1}-...-x_{i-1}\right)  -x_{i}=-z-x_{1}-...-x_{n-1}.
\]
In either case, the above condition is satisfied by $\left(  z,x,w^{\prime
}\right)  =\pi_{n,\partial F\ast}\left(  z,x,w\right)  \in\mathcal{F}%
^{n}|_{F_{n}}$, and hence $\pi_{n,\partial F\ast}\left(  z,x,w\right)
\in\widetilde{\frak{F}_{n-1}^{\prime}}$.

With $\widetilde{\frak{F}_{n-1}^{\prime}}$ canonically embedded in
$\mathcal{F}^{n}|_{F_{n}}$, the composition
\[
\eta_{n-1\ast}^{\prime\prime}:=\eta_{n-1\ast}^{\prime}\circ\pi_{n,\partial
F\ast}:\widetilde{\frak{F}_{n-1}^{\prime\prime}}\rightarrow\frak{F}%
_{n-1}^{\prime}
\]
of quotient maps realizes $\frak{F}_{n-1}^{\prime}$ as a quotient groupoid
$\widetilde{\frak{F}_{n-1}^{\prime\prime}}/\sim$ of $\widetilde{\frak{F}%
_{n-1}^{\prime\prime}}$ where the equivalence relation $\sim$ is the original
one extended to $\widetilde{\frak{F}_{n-1}^{\prime\prime}}$, i.e. the relation
generated by
\[
(z,x,w)\sim(z,x,w_{1},...,w_{i}=\infty,\infty,...,\infty)
\]
for all $(z,x,w)$ with $w_{i}=\infty$ for an $1\leq i\leq n$. Since
$\pi_{n,\partial F\ast}:\widetilde{\frak{F}_{n-1}^{\prime\prime}}%
\rightarrow\widetilde{\frak{F}_{n-1}^{\prime}}\subset\mathcal{F}^{n}|_{F_{n}}$
where the inclusion is implemented by $\iota_{n-1}^{\prime}:C^{\ast}\left(
\widetilde{\frak{F}_{n-1}^{\prime}}\right)  \rightarrow C^{\ast}\left(
\mathcal{F}^{n}|_{F_{n}}\right)  $ at the C*-algebra level, it is easy to see
that $\eta_{n-1}^{\prime\prime}=\pi_{n,\partial F}\circ\eta_{n-1}^{\prime}$
and hence
\[
C(S_{q}^{2n+1})|_{\partial F}=\pi_{n,\partial F}\left(  C(S_{q}^{2n+1}%
)|_{F_{n}}\right)
\]
\[
=\pi_{n,\partial F}\left(  \eta_{n-1}^{\prime}\left(  C^{\ast}(\frak{F}%
_{n-1}^{\prime})\right)  \right)  =\eta_{n-1}^{\prime\prime}\left(  C^{\ast
}(\frak{F}_{n-1}^{\prime})\right)  \subset C^{\ast}(\mathcal{F}^{n}|_{\partial
F}).
\]

Note that $\widetilde{\frak{F}_{n-1}^{\prime\prime}}\subset\mathcal{F}%
^{n}|_{\partial F}$ is a closed subgroupoid of $\widetilde{\frak{F}_{n}%
}\subset\mathcal{F}^{n}$ and hence $\frak{F}_{n-1}^{\prime}$ is a closed
subgroupoid of $\frak{F}_{n}$ with $\eta_{n-1\ast}^{\prime\prime}$ equal to
the restriction of $\eta_{n\ast}$ to $\widetilde{\frak{F}_{n-1}^{\prime\prime
}}$. On the other hand, on the `interior' $\widetilde{\frak{F}_{n}}%
-\widetilde{\frak{F}_{n-1}^{\prime\prime}}=\mathcal{F}^{n}-\left(
\mathcal{F}^{n}|_{\partial F}\right)  $, the relation $\sim$ is trivial and
hence
\[%
\begin{array}
[c]{ccccccccc}%
&  &  &  &  &  & C(S_{q}^{2n+1})|_{\partial F} &  & \\
&  &  &  &  &  & \parallel &  & \\
0 & \rightarrow & \iota_{n}\eta_{n}\left(  C^{\ast}(\frak{F}_{n}%
-\frak{F}_{n-1}^{\prime})\right)  & \subset & \iota_{n}\eta_{n}\left(
C^{\ast}(\frak{F}_{n})\right)  & \overset{\rho_{\partial F}}{\rightarrow} &
\eta_{n-1}^{\prime\prime}\left(  C^{\ast}(\frak{F}_{n-1}^{\prime})\right)  &
\rightarrow & 0\\
&  & \parallel &  & \cap &  & \cap &  & \\
0 & \rightarrow & \iota_{n}C^{\ast}(\widetilde{\frak{F}_{n}}-\widetilde
{\frak{F}_{n-1}^{\prime\prime}}) & \subset & \iota_{n}C^{\ast}(\widetilde
{\frak{F}_{n}}) & \overset{\rho_{\partial F}}{\rightarrow} & C^{\ast
}(\widetilde{\frak{F}_{n-1}^{\prime\prime}}) & \rightarrow & 0\\
&  & \parallel &  & \cap &  & \cap &  & \\
0 & \rightarrow &  C(\mathbb{T})\otimes\mathcal{K}\left(  \ell^{2}\left(
\mathbb{Z}_{\geq}^{n}\right)  \right)  & \subset &  C^{\ast}(\mathcal{F}%
^{n}) & \overset{\rho_{\partial F}}{\rightarrow} & C^{\ast}(\mathcal{F}%
^{n}|_{\partial F}) & \rightarrow & 0
\end{array}
\]

But as we have seen $C(S_{q}^{2n+1})$ is the only C*-subalgebra of
$C^{*}(\mathcal{F}^{n})$ whose image under $\rho_{\partial F}$ coincides with
$C(S_{q}^{2n+1})|_{\partial F}\subset C^{*}(\mathcal{F}^{n}|_{\partial F})$.
So we get $\iota_{n}\eta_{n}\left(  C^{*}(\frak{F}_{n})\right)  =C(S_{q}%
^{2n+1})\subset C^{*}(\mathcal{F}^{n})$ as stated.%

%TCIMACRO{\TeXButton{End Proof}{\endproof}}%
%BeginExpansion
\endproof
%EndExpansion%

%TCIMACRO{\TeXButton{BeginCorollary}{\begin{corollary}}}%
%BeginExpansion
\begin{corollary}%
%EndExpansion
The C*-algebra $C(S_{q}^{2n+1})$ is independent of $q$.%

%TCIMACRO{\TeXButton{EndCorollary}{\end{corollary}}}%
%BeginExpansion
\end{corollary}%
%EndExpansion

We remark that this fact can also be derived from Vaksman and Soibelman's
result \cite{VaSo}.

\end{document}